\documentclass[amstex,12pt,russian,amssymb]{article}

\usepackage{mathtext}
\usepackage[cp1251]{inputenc}
\usepackage[T2A]{fontenc}
\usepackage[russian]{babel}
\usepackage[dvips]{graphicx}
\usepackage{amsmath}
\usepackage{amssymb}
\usepackage{amsxtra}
\usepackage{latexsym}
\usepackage{ifthen}

\begin{document}

\newcounter{lemma}
\newcommand{\lemma}{\par \refstepcounter{lemma}%
{\bf Лемма \arabic{lemma}.}}

\newcounter{corollary}
\newcommand{\corollary}{\par \refstepcounter{corollary}%
{\bf Следствие \arabic{corollary}.}}

\newcounter{remark}
\newcommand{\remark}{\par \refstepcounter{remark}%
{\bf Замечание \arabic{remark}.}}

\newcounter{theorem}
\newcommand{\theorem}{\par \refstepcounter{theorem}%
{\bf Теорема \arabic{theorem}.}}

\newcounter{proposition}
\newcommand{\proposition}{\par \refstepcounter{proposition}%
{\bf Предложение \arabic{proposition}.}}

\renewcommand{\refname}{\centerline{\bf Список литературы}}

\newcommand{\proof}{{\it Доказательство.\,\,}}

\noindent УДК 517.5

{\bf Е.А.~Севостьянов} (Житомирский государственный университет
имени Ивана Франко)

{\bf Є.О.~Севостьянов} (Житомирський державний університет імені
Івана Фран\-ка)

{\bf E.A.~Sevost'yanov} (Zhytomyr Ivan Franko State University)

\medskip
{\bf Об изолированных особенностях отображений, обратные к которым
обобщённо квазиконформны}

{\bf Про ізольовані сингулярності відображень, обернені до яких
узагальнено квазіконформні}

{\bf On isolated singularities of mappings, inverse of which are
generalized  qua\-si\-con\-for\-mal}

\medskip\medskip
Установлено, что гомеоморфизмы областей евклидового пространства,
обратные к которым искажают модуль семейств кривых по типу
неравенства Полецкого, продолжаются по непрерывности в изолированную
точку границы. 

\medskip\medskip
Доведено, що гомеоморфізми областей евклідового простору, обернені
до яких спотворюють модуль сімей кривих по типу нерівності
Полецького, мають неперервне продовження в ізольовану точку межі. 

\medskip\medskip
We have proved that homeomorphisms of domains of Euclidean space,
inverse of which distort the modulus of families of curves by
Poletskii type, have a continuous extension to isolated boundary
point. 

\newpage
{\bf 1.~Введение.} В недавней работе~\cite{SevSkv} установлено, что
семейства гомеоморфизмов $g:D^{\,\prime}\rightarrow D,$ обратные к
которым удовлетворяют условию вида
\begin{equation}\label{eq1}
M(f(\Gamma))\leqslant \int\limits_DQ(x)\cdot
\rho^n(x)\,dm(x)\quad\forall\,\rho\in {\rm adm}\,\Gamma\,,
\end{equation}
$Q\in L^1(D),$ равностепенно непрерывны в $D^{\,\prime}$ (см.
\cite[теорема~2]{SevSkv}). Этот результат справедлив даже в общих
метрических пространствах, см. там же. В данной заметке мы
установим, что указанные гомеоморфизмы при тех же условиях имеют
непрерывное продолжение в изолированную точку границы области.
Доказательство этого результата мы приводим для некоторого более
общего класса отображений, нежели тех, что
удовлетворяют~(\ref{eq1}).

\medskip
Следует отметить работы, в которых данная проблема уже решалась в
менее полном объёме (см., напр., \cite[следствие~5.23]{MRSY},
\cite[теорема~6.1]{RSa}, \cite[теорема~5]{Sm}
и~\cite[теорема~6.1]{ARS}). Во всех этих работах случай
изолированной особой точки учитывался, поскольку изолированная точка
границы является слабо плоской~\cite[теорема~10.12]{Va}. С другой
стороны, область $D$ предполагалась локально связной на своей
границе, и это -- {\it лишнее предположение} для изолированных
точек. Ниже мы докажем, что гомеоморфизм $g$, являющийся обратным к
отображению $f:D\rightarrow D^{\,\prime}$ {\it произвольной} области
$D\subset {\Bbb R}^n$ на область $D^{\,\prime}\subset {\Bbb R}^n,$
$n\geqslant 2,$ продолжается по непрерывности в изолированную
граничную точку $x_0\in
\partial D^{\,\prime},$ если $f$ удовлетворяет~(\ref{eq1}), а
$Q$ интегрируема в $D.$ Никаких дополнительных условий на границу
области $D$ не требуется.

\medskip Пусть $x_0\in {\Bbb R}^n,$ $r>0,$ $0<r_1<r_2<\infty.$ Всюду далее, как обычно,
$$B(x_0,r)=\{ x\,\in\,{\Bbb R}^n : |x-x_0|<r\}\,,$$
$$S(x_0,r)=\{ x\,\in\,{\Bbb R}^n : |x-x_0|=r\}\,,\quad
A(x_0, r_1, r_2)=\{ x\,\in\,{\Bbb R}^n : r_1<|x-x_0|<r_2\}\,.$$
Пусть $Q:{\Bbb R}^n\rightarrow[0, \infty],$ $Q(x)\equiv 0$ при
$x\not\in D,$ -- измеримая по Лебегу функция. Гомеоморфизм
$f:D\rightarrow \overline{{\Bbb R}^n},$ $\overline{{\Bbb R}^n}={\Bbb
R}^n\cup\{\infty\},$ называется {\it кольцевым $Q$-гомеоморфизмом в
точке $x_0\in D,$} если существует $r_0=r_0(x_0)>0$ такое, что
соотношение
\begin{equation}\label{eq3*!}
M\left(f\left(\Gamma\left(S_1,\,S_2,\,A\right)\right)\right)\leqslant
\int\limits_{A} Q(x)\cdot \eta^n(|x-x_0|)\, dm(x) \end{equation}
выполнено для любого кольца $A=A(x_0, r_1, r_2),$ $0<r_1<r_2< r_0$ и
для каждой измеримой функции $\eta:(r_1,r_2)\rightarrow [0,\infty ]$
такой, что
\begin{equation}\label{eq*3}
\int\limits_{r_1}^{r_2}\eta(r)\, dr \geqslant 1\,.
\end{equation}
Свойство~(\ref{eq*3}) можно также определить в точке $x_0=\infty$
при помощи инверсии $\varphi=x/|x|^2.$
Гомеоморфизм $f:D\rightarrow \overline{{\Bbb R}^n}$ называется {\it
кольцевым $Q$-гомеоморфизмом в $E\subset \overline{D},$} если
соотношение~(\ref{eq3*!}) выполнено в каждой точке $x_0\in E.$
Основной результат настоящей заметки можно сформулировать следующим
образом.

\medskip
\begin{theorem}\label{th1}
{\sl\, Пусть $D, D^{\,\prime}$ -- области в $\overline{{\Bbb R}^n},$
$n\geqslant 2,$ и пусть $g:D^{\,\prime}\rightarrow D,$
$g(D^{\,\prime})=D$ -- гомеоморфизм, обратный $f:=g^{\,-1}$ к
которому удовлетворяет условию ~(\ref{eq3*!}) в каждой точке $x_0\in
\partial D.$ Если $Q\in L^1(D)$ и $y_0$ -- изолированная точка границы
области $D^{\,\prime},$ то отображение $g$ имеет непрерывное
продолжение $\overline{g}:D^{\,\prime}\cup\{y_0\}\rightarrow
\overline{{\Bbb R}^n}$ в точку $y_0.$ }
\end{theorem}

\medskip
Здесь и далее граница, замыкание и непрерывное продолжение
отображений понимаются в смысле расширенного
пространства~$\overline{{\Bbb R}^n}.$

\medskip
{\bf 2.~Доказательство основного результата.} Для отображения
$f:D\rightarrow \overline{{\Bbb R}^n}$ и множества $E\subset
\overline{D}\subset\overline{{\Bbb R}^n}$ положим
$$C(f, E)=\{y\in \overline{{\Bbb R}^n}:\exists\, x\in E, x_k\in D: x_k\rightarrow x, f(x_k)
\rightarrow y, k\rightarrow\infty\}\,.$$ Следующее важнейшее
утверждение опубликовано в~\cite[лемма~5.3]{IR$_1$}, см.
также~\cite[лемма~6.5]{MRSY$_1$}.

\medskip
\begin{proposition}\label{pr1}
{\sl\, Пусть $D$ -- область в $\overline{{\Bbb R}^n},$ $n\geqslant
2,$ и пусть $f:D\rightarrow \overline{{\Bbb R}^n}$ -- гомеоморфизм.
Тогда существует взаимно однозначное соответствие между компонентами
$K$ и $K^{\,\prime}$ границ $\partial D$ и $\partial D^{\,\prime}$
такое, что $C(f, K)=K^{\,\prime}$ и $C(f^{\,-1}, K^{\,\prime})=K.$}
\end{proposition}

\medskip
Следующая лемма указывает на выполнение условия <<слабой плоскости>>
в изолированных граничных точках произвольной области.

\medskip
\begin{lemma}\label{lem2}
Пусть $D$ -- область в ${\Bbb R}^n,$ $n\geqslant 2,$ и $x_0\in D.$
Тогда для каждого $P>0$ и для любой окрестности $U$ точки $x_0$
найдётся окрестность $V\subset U$ этой же точки такая, что
$M(\Gamma(E, F, D\setminus\{x_0\}))>P$ для произвольных континуумов
$E, F\subset D\setminus\{x_0\},$ пересекающих $\partial U$ и
$\partial V.$
\end{lemma}

\medskip
\begin{proof}
Пусть $U$ -- произвольная окрестность точки $x_0.$ Выберем
$\varepsilon_0>0$ так, чтобы $\overline{B(x_0,
\varepsilon_0)}\subset D\cap U.$ Пусть $c_n$ -- положительная
постоянная, определённая в соотношении~(10.11) в \cite{Va}, а число
$\varepsilon\in(0, \varepsilon_0)$ настолько мало, что
$c_n\cdot\log\frac{\varepsilon_0}{\varepsilon}>P.$ Положим
$V:=B(x_0, \varepsilon).$ Пусть $E, F$ -- произвольные континуумы,
пересекающие $\partial U$ и $\partial V,$ тогда также $E$ и $F$
пересекают $S(x_0, \varepsilon_0)$ и $\partial V$ (см.
\cite[теорема~1.I, гл.~5, \S\, 46]{Ku}). Необходимое заключение
вытекает на основании~\cite[разд.~10.12]{Va}, поскольку
$M(\Gamma(E, F, D\setminus\{x_0\}))\geqslant
c_n\cdot\log\frac{\varepsilon_0}{\varepsilon}>P.\quad \Box$
\end{proof}

\bigskip
{\it Доказательство теоремы~\ref{th1}.} Без ограничения общности
можно считать, что $y_0\ne\infty.$ Обозначим далее через $h(x, y)$
хордальное (сферическое) расстояние между точками $x, y\in
\overline{{\Bbb R}^n},$
$$h(x,\infty)=\frac{1}{\sqrt{1+{|x|}^2}}\,,\quad
h(x,y)=\frac{|x-y|}{\sqrt{1+{|x|}^2} \sqrt{1+{|y|}^2}}\,,\quad  x\ne
\infty\ne y\,.$$
Предположим противное, а именно, предположим, что $g$ не имеет
предела в точке $y_0.$ Поскольку $\overline{{\Bbb R}^n}$ -- компакт,
$C(g, y_0)\ne\varnothing.$ Тогда найдутся $x_1, x_2\in
\overline{{\Bbb R}^n},$ $x_1\ne x_2,$ и не менее двух
последовательностей $y_m, y^{\,\prime}_m\rightarrow y_0,$
$m\rightarrow\infty,$ таких что $z_m:=g(y_m)\rightarrow x_1,$
$z_m^{\,\prime}=g(y^{\,\prime}_m)\rightarrow x_2$ при
$m\rightarrow\infty.$ 
Не ограничивая общности, можно считать, что $x_1\ne\infty.$

\medskip
По предложению~\ref{pr1} $C(g, y_0)$ -- континуум в $\overline{{\Bbb
R}^n}$ и является компонентой границы~$\partial D.$ Покажем, что
найдётся $\varepsilon_1>0$ такое, что
\begin{equation}\label{eq2}
B(x_1, \varepsilon_1)\cap K=\varnothing
\end{equation}
для всякой компоненты $K$ границы $\partial D$ такой, что $C(g,
y_0)\ne K.$ Предположим противное. Тогда найдётся последовательность
компонент границы $d_m\subset
\partial D,$ $d_m\ne C(g, y_0),$ $m=1,2,\ldots ,$ такая, что $B(x_1, 1/m)\cap d_m\ne
\varnothing.$ По предложению~\ref{pr1} компоненте $d_m$ границы
области $D$ соответствует компонента $d^{\,\prime}_m\subset
\partial D^{\,\prime}$ так, что $C(f, d_m)=d^{\,\prime}_m.$
Поэтому можно выбрать $\zeta_m\in B(x_1, 1/m)\cap D$ так, что
$$h(f(\zeta_m),  d^{\,\prime}_m)=\inf\limits_{p\in
d^{\,\prime}_m}h(f(\zeta_m), p)<1/m\,.$$
Так как $d^{\,\prime}_m$ -- компакт в $\overline{{\Bbb R}^n},$
найдётся $\xi_m\in d^{\,\prime}_m$ такое, что $h(f(\zeta_m),
d^{\,\prime}_m)=h(f(\zeta_m), \xi_m).$ Поскольку $\overline{{\Bbb
R}^n}$ -- компакт, то существует подпоследовательность
$f(\zeta_{m_k}),$ сходящаяся в $\overline{{\Bbb R}^n}$ при
$k\rightarrow\infty.$ Поскольку $\zeta_{m_k}\in B(x_1, 1/m_k),$ то
по предложению~\ref{pr1} последовательность $f(\zeta_{m_k})$ может
сходиться только к $y_0$ при $k\rightarrow\infty.$ Тогда по
неравенству треугольника
$$h(y_0, d^{\,\prime}_{m_k})\leqslant h(y_0, \xi_{m_k})\leqslant
h(y_0, f(\zeta_{m_k}))+ h(f(\zeta_{m_k}),
\xi_{m_k})\stackrel{k\rightarrow\infty}\rightarrow 0\,,$$
что противоречит изолированности особенности $y_0.$ Полученное
противоречие указывает на справедливость формулы~(\ref{eq2}).

\begin{figure}[h]
\centerline{\includegraphics[scale=0.5]{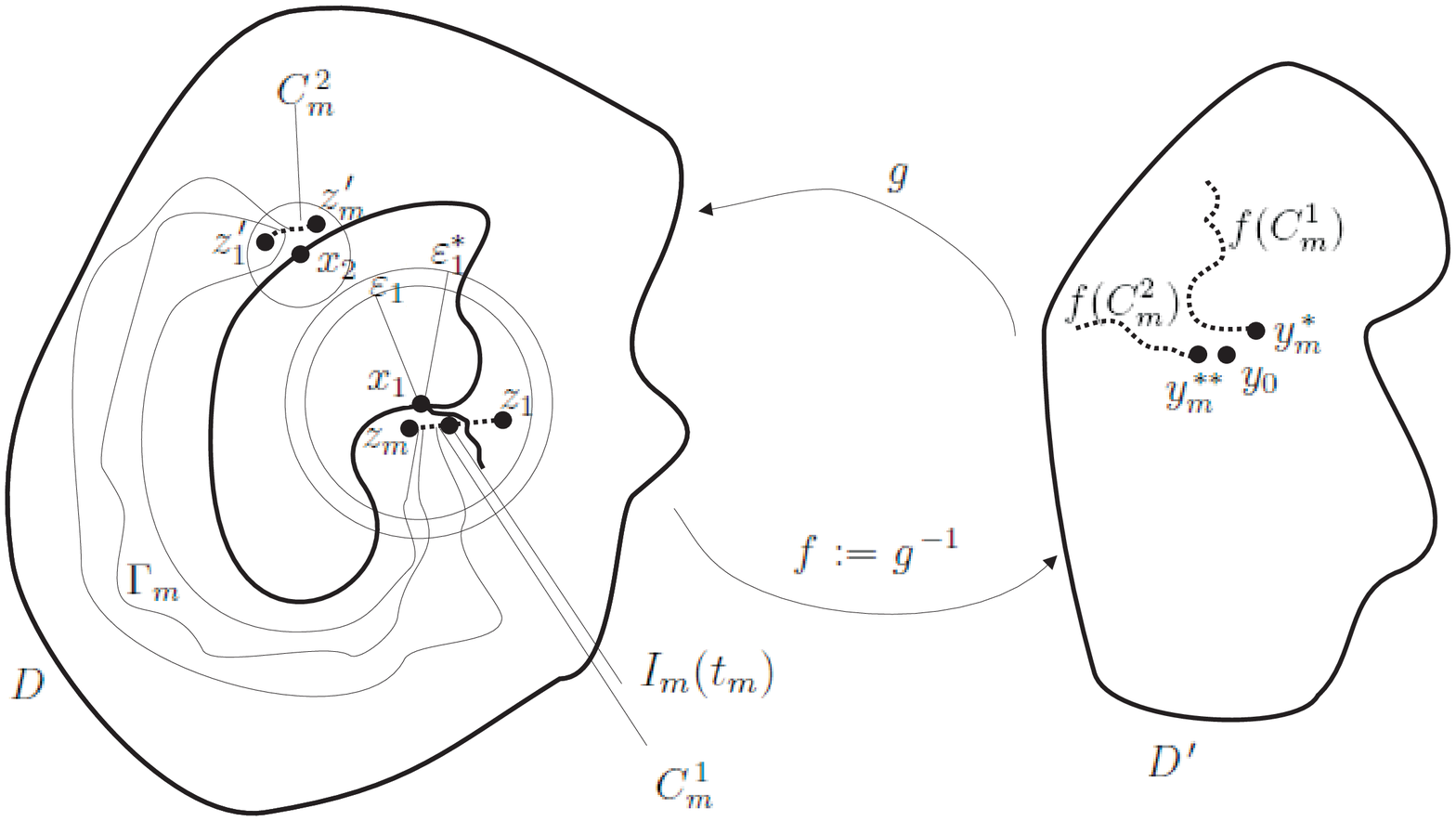}} \caption{ К
доказательству теоремы~\ref{th1}}\label{fig1}
\end{figure}

\medskip
Пусть $B_*(x_2, \varepsilon_2)=B(x_2, \varepsilon_2)$ при $x_2\ne
\infty$ и $B_*(x_2, \varepsilon_2)=\{x\in \overline{{\Bbb R}^n}:
h(x, \infty)<\varepsilon_2\}$ при $x_2=\infty.$ Ввиду соображений,
аналогичных тем, что изложены при доказательстве
соотношения~(\ref{eq2}), найдётся $\varepsilon_2>0$ такое, что
$$B_*(x_2, \varepsilon_2)\cap K=\varnothing$$
для всякой компоненты $K$ границы $\partial D,$ не равной $C(g,
y_0)\ne K.$
Без ограничения общности, мы можем считать, что $\overline{B(x_1,
\varepsilon_1)}\cap \overline{B_*(x_2, \varepsilon_2)}=\varnothing,$
$z_m\in B(x_1, \varepsilon_1)$ и $z^{\,\prime}_m \in B_*(x_2,
\varepsilon_2)$ при всех $m=1,2,\ldots $ (см. рисунок~\ref{fig1}).
Заметим, что множество~$B(x_1, \varepsilon_1)$ выпукло, а $B_*(x_2,
\varepsilon_2)$ линейно связно. В таком случае точки $z_1$ и $z_m$
можно соединить отрезком $I_m(t)=z_1+t(z_m-z_1),$ $t\in (0, 1),$
лежащим внутри шара~$B(x_1, \varepsilon_1).$ Аналогично, точки
$z^{\,\prime}_1$ и $z^{\,\prime}_m$ можно соединить кривой
$J_m=J_m(t),$ $t\in [0, 1],$ лежащей в <<шаре>> $B_*(x_2,
\varepsilon_2).$

\medskip
Всюду далее для кривой $\alpha:[a, b]\rightarrow {\Bbb R}^n$ мы
используем обозначение
$$|\alpha|:=\{x\in {\Bbb R}^n: \exists\,t\in
[a, b]: \alpha(t)=x\}\,.$$
Заметим, что множество $|I_m|$ не обязано лежать в области $D$
целиком. Тем не менее, в этом случае, найдётся $t_m\in [0, 1]$
такое, что $h(I_m(t_m),
\partial D)<1/m$ и $|I_m|_{[0, t_m]}|\subset D.$ Аналогично, $J_m$ может не лежать
в области $D$ целиком, однако, найдётся $p_m\in [0, 1]$ такое, что
$h(J_m(p_m),
\partial D)<1/m$ и $|J_m|_{[0, p_m]}|\subset D.$
Если $I_m\subset D$ либо $J_m\subset D,$ то полагаем $t_m:=1$ и
$p_m:=1,$ соответственно. Положим $C^{\,1}_m:=I_m|_{[0, t_m]}$ и
$C^{\,2}_m:=J_m|_{[0, p_m]}.$
Рассмотрим последовательности~$y^{\,*}_m:=f(I_m(t_m)),$
$y^{\,**}_m:~=~f(J_m(p_m)).$
Поскольку пространство~$\overline{{\Bbb R}^n}$ компактно, мы можем
считать, что все рассматриваемые последовательности $I_m(t_m),$
$J_m(p_m),$ $y^{\,*}_m$ и $y^{\,**}_m$ сходятся при
$m\rightarrow\infty.$

\medskip
Покажем, что $y^{\,*}_m, y^{\,**}_m\rightarrow y_0$ при
$m\rightarrow\infty.$ Пусть $y^{\,*}_m\rightarrow w_0$ при
$m\rightarrow\infty.$ Поскольку последовательность $I_m(t_m)$
сходится по предположению и, при этом, $h(I_m(t_m),
\partial D)<1/m,$ то $I_m(t_m)$ сходится к некоторой точке $\omega_0\in\partial D.$
Так как $I_m(t_m)\in B(x_1, \varepsilon_1),$ то $\omega_0\in C(g,
y_0),$ поскольку в шаре $B(x_1, \varepsilon_1)$ нет других компонент
границы $\partial D,$ отличных от $C(g, y_0),$ ввиду
соотношения~(\ref{eq2}). Тогда из условий $y^{\,*}_m=f(I_m(t_m))$ и
$y^{\,*}_m\rightarrow y_0,$ $m\rightarrow\infty,$ вытекает, что
$w_0\in C(f, C(g, y_0)).$ Так как $g=f^{\,-1},$ то по
предложению~\ref{pr1} $w_0=y_0,$ что и требовалось установить.
Рассуждая аналогично, мы можем показать, что $y^{\,**}_m\rightarrow
y_0$ при $m\rightarrow\infty.$

\medskip
Поскольку $\overline{B(x_1, \varepsilon_1)}\cap \overline{B_*(x_2,
\varepsilon_2)}=\varnothing,$ то при некотором
$\varepsilon^*_1>\varepsilon_1$ мы также всё ещё имеем
$\overline{B(x_1, \varepsilon^*_1)}\cap \overline{B_*(x_2,
\varepsilon_2)}=\varnothing.$ Пусть $\Gamma_m=\Gamma(|C^1_m|,
|C^2_m|, D).$ Заметим, что
\begin{equation}\label{eq3}
\Gamma_m>\Gamma(S(x_1, \varepsilon^*_1), S(x_1, \varepsilon_1),
A(x_1, \varepsilon_1, \varepsilon^*_1))\,.
\end{equation}
В самом деле, пусть $\gamma\in \Gamma_m,$ $\gamma:[a, b]\rightarrow
{\Bbb R}^n.$ Поскольку $\gamma(a)\in |C^1_m|\subset B(x_0,
\varepsilon_1)$ и $\gamma(b)\in |C^2_m|\in {\Bbb R}^n\setminus
B(x_0, \varepsilon_1),$ то ввиду~\cite[теорема~1.I, гл.~5, \S\,
46]{Ku} найдётся $t_1\in (a, b)$ такое, что $\gamma(t_1)\in S(x_1,
\varepsilon_1).$ Не ограничивая общности, можно считать, что
$|\gamma(t)-x_1|>\varepsilon_1$ при всех $t>t_1.$ Далее, так как
$\gamma(t_1)\in B(x_1, \varepsilon^*_1)$ и $\gamma(b)\in
|C^2_m|\subset {\Bbb R}^n\setminus B(x_0, \varepsilon^*_1),$ то
ввиду~\cite[теорема~1.I, гл.~5, \S\, 46]{Ku} найдётся $t_2\in (t_1,
b)$ такое, что $\gamma(t_2)\in S(x_1, \varepsilon^*_1).$ Не
ограничивая общности можно считать, что
$|\gamma(t)-x_1|<\varepsilon_1^*$ при $t_1<t<t_2.$ Таким образом,
$\gamma|_{[t_1, t_2]}$ -- подкривая кривой $\gamma,$ принадлежащая
семейству~$\Gamma(S(x_1, \varepsilon^*_1), S(x_1, \varepsilon_1),
A(x_1, \varepsilon_1, \varepsilon^*_1)).$ Таким образом,
соотношение~(\ref{eq3}) установлено.
Рассмотрим функцию
$$\eta(t)\quad =\quad\left\{
\begin{array}{rr}
1/(\varepsilon^{*}_1-\varepsilon_1), & t\in [\varepsilon_1, \varepsilon^*_1],\\
0, & t\in {\Bbb R}\setminus[\varepsilon_1,
\varepsilon^*_1]\,.\end{array} \right.$$
Заметим, что функция $\eta$ удовлетворяет соотношению~(\ref{eq*3})
при $r_1=\varepsilon_1$ и $r_2=\varepsilon^*_1.$ Тогда по
определению кольцевого $Q$-отображения в точке $x_1,$ с учётом
условия $Q\in L^1(D)$ и по соотношению~(\ref{eq3}), получаем:
\begin{equation}\label{eq4}
M(f(\Gamma_m))\leqslant M(f(\Gamma(S(x_1, \varepsilon^*_1), S(x_1,
\varepsilon_1), A(x_1, \varepsilon_1, \varepsilon^*_1))))\leqslant
\Vert Q\Vert/(\varepsilon^{*}_1-\varepsilon_1)^n <\infty\,,
\end{equation}
где $\Vert Q\Vert$ обозначает $L^1$-норму функции $Q$ в области $D.$
Для множеств $A, B\subset{\Bbb R}^n$ положим, как обычно,
$${\rm diam}\,A=\sup\limits_{x, y\in A}|x-y|\,,\quad {\rm dist}\,(A, B)=\inf\limits_{x\in A,
y\in B}|x-y|\,.$$
Покажем, что соотношение~(\ref{eq4}) противоречит условию слабой
плоскости в точке $y_0$ (см. лемму~\ref{lem2}). В самом деле,
$${\rm diam}\,f(C^1_m)\geqslant |f(z_1)-f(I_m(t_m))|=|y_1-y^*_m| \geqslant
(1/2)\cdot|y_1-y_0|>0$$
и
$${\rm diam}\,f(C^2_m)\geqslant |f(z^{\,\prime}_1)-f(J_m(p_m))|
=|y_1^{\,\prime}-y^{**}_m| \geqslant
(1/2)\cdot|y^{\,\prime}_1-y_0|>0$$
при больших $m\in {\Bbb N},$ кроме того,
$${\rm dist}\,(f(C^1_m), f(C^2_m))\leqslant |y^*_m-y^{**}_m|\rightarrow 0,\quad m\rightarrow
\infty\,.$$ Тогда ввиду леммы~\ref{lem2}
$$M(f(\Gamma_m))=M(\Gamma(f(C^1_m), f(C^2_m), D^{\,\prime}))
\rightarrow\infty\,,\quad m\rightarrow\infty\,,$$
что противоречит соотношению~(\ref{eq4}). Полученное противоречие
указывает на ошибочность исходного предположения о том, что
отображение $g$ не имеет непрерывного продолжения в
точку~$y_0.$~$\Box$

\medskip
{\bf Пример.} Пусть $p\geqslant 1,$ $\alpha\in \left(0,
n/p(n-1)\right)$ и $e_1=(0,0,\ldots, 0, 1/2).$ Зададим отображение
$f$ в области $D:={\Bbb B}^n\setminus\{e_1\cup 0\}$ следующим
образом: $$f(x)~=~\frac{1+|x|^{\alpha}}{|x|}\cdot x\,,\quad x\in
{\Bbb B}^n\setminus\{e_1\cup 0\}\,.$$
Нетрудно видеть, что отображение $f$ является кольцевым
$Q$-гомеоморфизмом ${\Bbb B}^n\setminus\{e_1\cup 0\}$ на
$A:=\{1<|y|<2\}\setminus\{e_2\},$ где
$Q(x):=\left(\frac{1+r^{\,\alpha}}{\alpha
r^{\,\alpha}}\right)^{n-1},$ $r=|x|$ (см., напр.,
\cite[предложение~6.3]{MRSY}). Более того,  $Q\in L^p({\Bbb B}^n).$
Заметим, что $f(e_1)=(0,0,\ldots,3/2):=e_2.$ Обратное отображение
$g:=f^{\,-1}(y)=\frac{y}{|y|}(|y|-1)^{1/\alpha}$ имеет непрерывное
продолжение в точку $e_2,$
$\overline{g}:=\frac{y}{|y|}(|y|-1)^{1/\alpha},$
$\overline{g}:\{1<|y|<2\}\rightarrow {\Bbb B}^n\setminus\{0\}$
(наличие этого продолжения также следует из теоремы~\ref{th1}). С
другой стороны, отображение $\overline{f}:=\overline{g^{\,-1}},$
$\overline{f}:{\Bbb B}^n\setminus\{0\}\rightarrow \{1<|y|<2\},$ не
имеет непрерывного продолжения в точку $0,$ которая является
изолированной точкой границы области ${\Bbb B}^n\setminus\{0\}.$
Последнее обстоятельство связано с неинтегрируемостью некоторой
функции $Q^{\,*}(y),$ отвечающей отображению $\overline{g}$ в
области $\{1<|y|<2\}$ в контексте неравенства~(\ref{eq3*!}).

\medskip
\begin{remark}\label{rem1}
Нетрудно видеть, что утверждение теоремы~\ref{th1} верно при
значительно более мягком условии на отображение $f,$ а именно,
достаточно потребовать условие~(\ref{eq3*!}) всего лишь в одной
конечной точке предельного множества~$C(g, y_0).$

\medskip
Более того, не очень существенно, требовать ли условие~(\ref{eq3*!})
на $\partial D$ или в $D.$ В самом деле, пусть в условиях
теоремы~\ref{th1} мы требуем соотношение~(\ref{eq3*!}) не на
$\partial D,$ а в каждой внутренней точке $x_0\in D.$ Повторяя
доказательство этой теоремы в тех же обозначениях, мы приходим к
cоотношению~(\ref{eq3}). Пусть теперь $a_k\in D,$ $k=1,2,\ldots ,$
некоторая (произвольная) последовательность точек, сходящаяся к
$x_1$ при $k\rightarrow\infty,$ такая что $|a_k-x_1|<1/k.$
Зафиксируем $x\in B(x_1, \varepsilon_1).$  Тогда по неравенству
треугольника $|x-a_k|\leqslant |x-x_1|+|x_1-a_k|<\varepsilon_1+1/k$
и, значит, $B(x_1, \varepsilon_1)\subset B(a_k, \varepsilon_1+1/k).$
Далее, для $x\in B(a_k, \varepsilon_1+2/k)$ по неравенству
треугольника имеем
$$|x-x_1|\leqslant |x-a_k|+|a_k-x_1|<\varepsilon_1+3/k\,.$$
Пусть $k_0\in {\Bbb N}$ настолько велико, что
$\varepsilon_1+3/k<\varepsilon_1^{\,*}$ при $k>k_0.$ Тогда $B(a_k,
\varepsilon_1+2/k)\subset B(x_1, \varepsilon_1^{*})$ при $k>k_0.$
Полагая $\widetilde{\varepsilon_1}:=\varepsilon_1+1/(k_0+1)$ и
$\widetilde{\varepsilon_2}:=\varepsilon_1+2/(k_0+1),$ получаем:
\begin{equation}\label{eq5}
B(x_1, \varepsilon_1)\subset B(a_{k_0+1},
\widetilde{\varepsilon_1})\subset B(a_{k_0+1},
\widetilde{\varepsilon_2})\subset B(x_1, \varepsilon_1^{\,*})\,.
\end{equation}
Рассуждая по аналогии с доказательством формулы~(\ref{eq3}), из
соотношения~(\ref{eq5}) получаем, что
\begin{equation}\label{eq6}\Gamma(S(x_1, \varepsilon^*_1), S(x_1, \varepsilon_1), A(x_1,
\varepsilon_1, \varepsilon^*_1))
>\Gamma(S(a_{k_0+1}, \widetilde{\varepsilon_1}), S(a_{k_0+1}, \widetilde{\varepsilon_2}),
A(a_{k_0+1}, \widetilde{\varepsilon_1},
\widetilde{\varepsilon_2}))\,,
\end{equation}
см. рисунок~\ref{fig2}.
\begin{figure}[h]
\centerline{\includegraphics[scale=0.5]{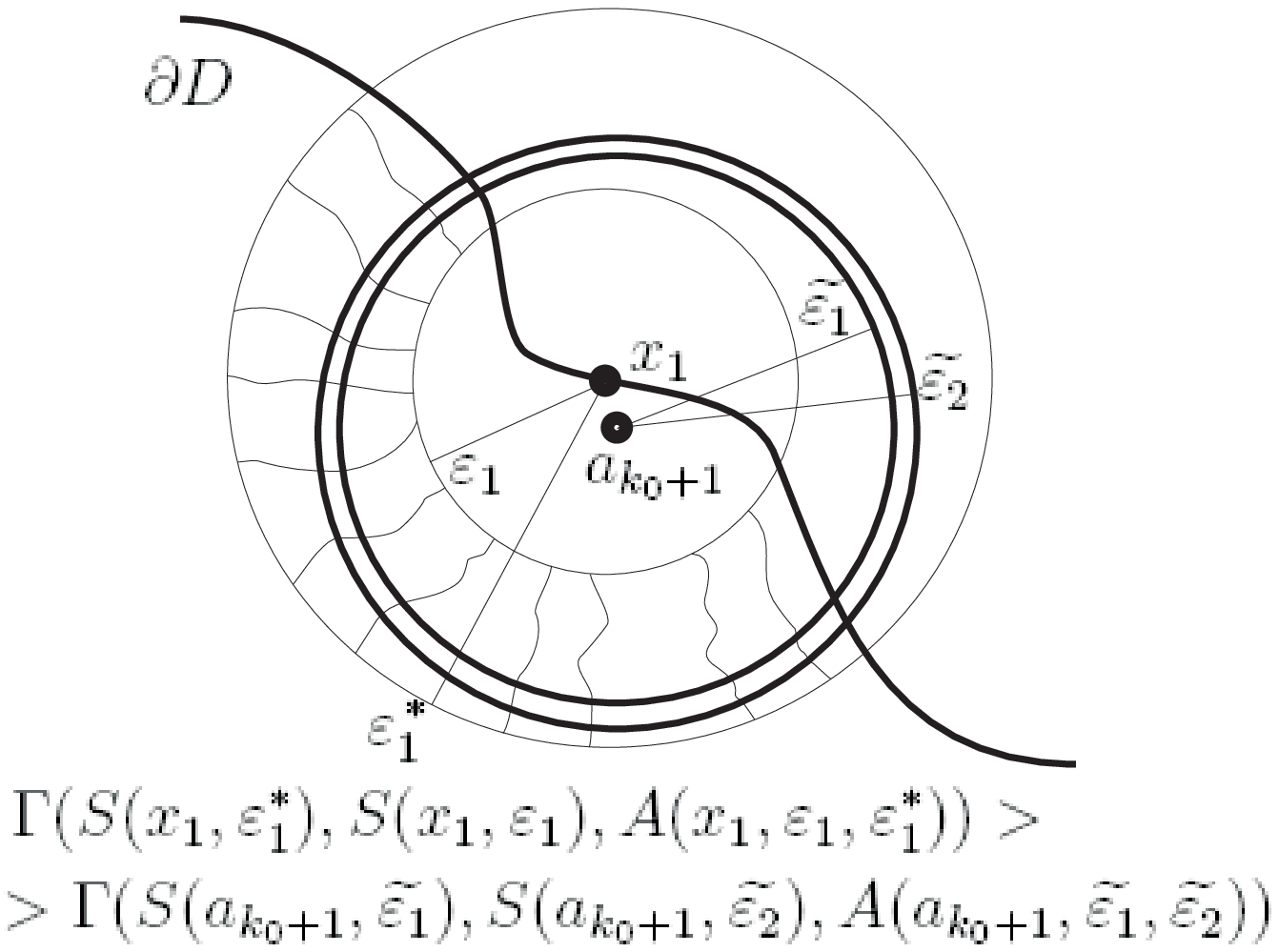}} \caption{К
замечанию~\ref{rem1}}\label{fig2}
\end{figure}
Рассмотрим функцию
$$\eta(t)\quad =\quad\left\{
\begin{array}{rr}
1/(\widetilde{\varepsilon_2}-\widetilde{\varepsilon_1}), & t\in [\widetilde{\varepsilon_1},
\widetilde{\varepsilon_2}],\\
0, & t\in {\Bbb R}\setminus[\widetilde{\varepsilon_1},
\widetilde{\varepsilon_2}]\,.\end{array} \right.$$
Заметим, что функция $\eta$ удовлетворяет соотношению~(\ref{eq*3})
при $r_1=\widetilde{\varepsilon_1}$ и
$r_2=\widetilde{\varepsilon_2}.$ Тогда по определению кольцевого
$Q$-отображения в точке $a_{k_0+1},$ с учётом условия $Q\in L^1(D)$
и по соотношениям~(\ref{eq3}) и~(\ref{eq6}), получаем:
$$M(f(\Gamma_m))\leqslant M(f(\Gamma(S(a_{k_0+1},
\widetilde{\varepsilon_1}), S(a_{k_0+1}, \widetilde{\varepsilon_2}),
A(a_{k_0+1}, \widetilde{\varepsilon_1},
\widetilde{\varepsilon_2}))))\leqslant$$
\begin{equation}\label{eq4*}
\leqslant \Vert
Q\Vert/(\widetilde{\varepsilon_2}-\widetilde{\varepsilon_1})^n
<\infty\,,
\end{equation}
где $\Vert Q\Vert$ обозначает $L^1$-норму функции $Q$ в области $D.$
Итак, вместо соотношения~(\ref{eq6}) мы имеем
соотношение~(\ref{eq4*}). Оставшаяся часть доказательства,
основанная на противоречии~(\ref{eq4*}) с условием слабой плоскости
в точке $y_0,$ не изменится.
\end{remark}

КОНТАКТНАЯ ИНФОРМАЦИЯ

\medskip
\noindent{{\bf Евгений Александрович Севостьянов} \\
Житомирский государственный университет им.\ И.~Франко\\
кафедра математического анализа, ул. Большая Бердичевская, 40 \\
г.~Житомир, Украина, 10 008 \\ тел. +38 066 959 50 34 (моб.),
e-mail: esevostyanov2009@gmail.com}

\end{document}